\documentclass[11pt,a4paper]{amsart}

\usepackage{tikz}
\usetikzlibrary{decorations.pathmorphing}
\usetikzlibrary{shapes.geometric}
\usepackage{fontawesome}
\usepackage{caption}
\usepackage{subcaption}

\usepackage[utf8]{inputenc}
\usepackage{amsmath,amssymb}

\usepackage[USenglish]{babel}

\begin{document}

\title[Computational Brainsphere Model Alzheimer's Disease Simulation]{A Computational Brainsphere Model for the Simulation of Alzheimer's Disease}

\author{Angela Kunoth}
\address{University of Cologne, Germany}
\email{akunoth@uni-koeln.de}
\thanks{The presented model was developed in joint work with Prof.\ Dr.\ Alexander Drzezga
    and Prof.\ Dr.\ Yaping Shao as part of the project ``Neurodegeneration Forecasting --- A Computational
    Brainsphere Model for the Simulation of Alzheimer's Disease'', supported by the Excellence Initiative of
    the University of Cologne, Germany, 2017-2019.}
\author{Anna Weller}
\thanks{Supported by the Excellence Initiative of the University of Cologne and by
  HESCOR (Human and Earth System Coupled Research) at the University of Cologne,
  part of the ``Profile Building 2022'' initiative by North Rhine-Westphalia's Ministry for Culture and Science,
Germany.}
\author{Tolunay Yilmaz}

\begin{abstract}
  Nowadays, there are different medical imaging techniques to collect data about the progression of
  Alzhei\-mer's disease in a patient's brain. These data describe different phenomena which are still
  not understood from a biological point of view.
  How can these data sets be combined in a mathematical model to simulate the evolution of such a
  neurodegenerative disease in a computer?
  In this snapshot, we present one possible approach to address this task with the help of graph theory
  and partial differential equations.
\end{abstract}

\maketitle

\section{Introduction}
Alzheimer's disease (AD) is a slowly progressing neurodegenerative disease causing the brain to shrink and brain cells to die.
Much progress has been made in investigating the course and certain biomarkers of the
disease. The perhaps largest attention is currently focused on two proteins that behave conspicuously
in the brain of an Alzheimer's patient: On the one hand, beta-amyloid peptides
aggregate in the surrounding area of the neuronal network in the brain as plaques \cite{EJ}. On the other
hand, there are abnormally folded tau proteins which accumulate in the neurons in form of \textit{tangles} and are
believed to affect their functioning and, finally, to cause their death \cite{GEC}. 

In the past, it was only possible to examine the brains of dead Alzheimer's patients.
The development of modern medical imaging techniques made it possible to visualize
protein plaques in the brains of living patients as well (in-vivo data). 
PET scanners measure accumulations of tau and beta-amyloid, and
functional MRI as well as DTI scanners provide indications of the structure and functioning of a
human brain. 
Central goal of subsequent studies is the behavior of tau and beta-amyloid as well as
the connectivity patterns in the brain.

A relatively new model for tau propagation is the transneuronal spread hypothesis. 
Transneuronal means that tangles travel from neuron to neuron, inducing tangles in neighboring neurons in a prion-like fashion. 
The word prion is a combination of  ``protein'', and ``infection'', indicating that missfolded proteins can infect others.
According to this hypothesis, the connectivity of neurons is a driver for tau pathology, i.e., highly interconnected neurons are more likely to receive pathology \cite{CR}.

On the other hand, it is also believed that aggregated tau proteins affect the functioning of neurons and their connections to other neurons in the brain. 
Additionally, the presence of beta-amyloid plaques in the surrounding of a neuron could influence its tendency to aggregate tau. 

Having discussed only a fraction of the possible interactions and factors of the disease process,
we already realize that it quickly becomes complicated to adequately describe such a
complex interplay. The crucial challenge is to analyze and combine the in-vivo data
within in a mathematical framework. We propose a \textit{Computational Brainsphere Model for
  the Simulation of Alzheimer's Disease} which consists of three interacting parts:
\begin{enumerate}
\item[1.] A \emph{Brain Network Model} to represent the human brain as a dynamical network, or,
  in mathematical terms, as a \emph{graph}, and to investigate the effects of tau pathology on the connectivity patterns;
\item[2.] A \emph{Tau Distribution Model} to simulate the spreading of tau proteins across the brain network;
\item[3.] A \emph{Beta-Amyloid Model} to model the dispersion of beta-amyloid in the environment of the brain network.
\end{enumerate}
In the following, we will give a brief introduction to the Brain Network and the Tau Distribution Model.
These have, as described above, a reciprocal dependency, i.e., the network structure influences the propagation of tau,
and the presence of tau influences the structure of the network.

\section{Graph Theory}\label{sec:graph-theory}
One major tool in the analysis of networks is a mathematical branch called \textit{graph theory}. It provides a formalization of a network as well as characteristics and algorithms to investigate structures of networks. 

\begin{figure}[ht]
\begin{center}
\begin{tikzpicture}
    \node[shape=circle,draw=black] (A) at (0,0) {$v_1$};
    \node[shape=circle,draw=black] (B) at (4,0) {$v_2$};
    \node[shape=circle,draw=black] (C) at (-1,1.5) {$v_3$};
    \node[shape=circle,draw=black] (D) at (-1.5,-2) {$v_4$};
    \node[shape=circle,draw=black] (E) at (-4,-0.5) {$v_5$};
    \path [-] (A) edge node[below] {$e_1$} (B);
    \path [-](A) edge node[left] {$e_2$} (C);
    \path [-](A) edge node[right] {$e_3$} (D);
    \path [-](B) edge node[above] {$e_4$} (C);
    \path [-](E) edge node[above] {$e_5$} (C);
    \path [-](D) edge node[below] {$e_6$} (E);
 \node[align=left] at (3.5,-1.5) {\footnotesize vertex $v_1$ is adjacent to $v_2, v_3$ and $v_4$ \\ \footnotesize vertex $v_1$ has degree $d_{v_1}=3$\\ \footnotesize edge $e_1$ is incident to $v_1$ and $v_2$};
\end{tikzpicture}
\caption{Example of a graph $G=(V,E)$ with vertex set $V= \lbrace v_1, v_2, v_3, v_4, v_5 \rbrace$ and edge set $E= \lbrace e_1, e_2, e_3, e_4, e_5, e_6 \rbrace$.} \label{fig:ex-graph}
 \end{center}
 \end{figure}
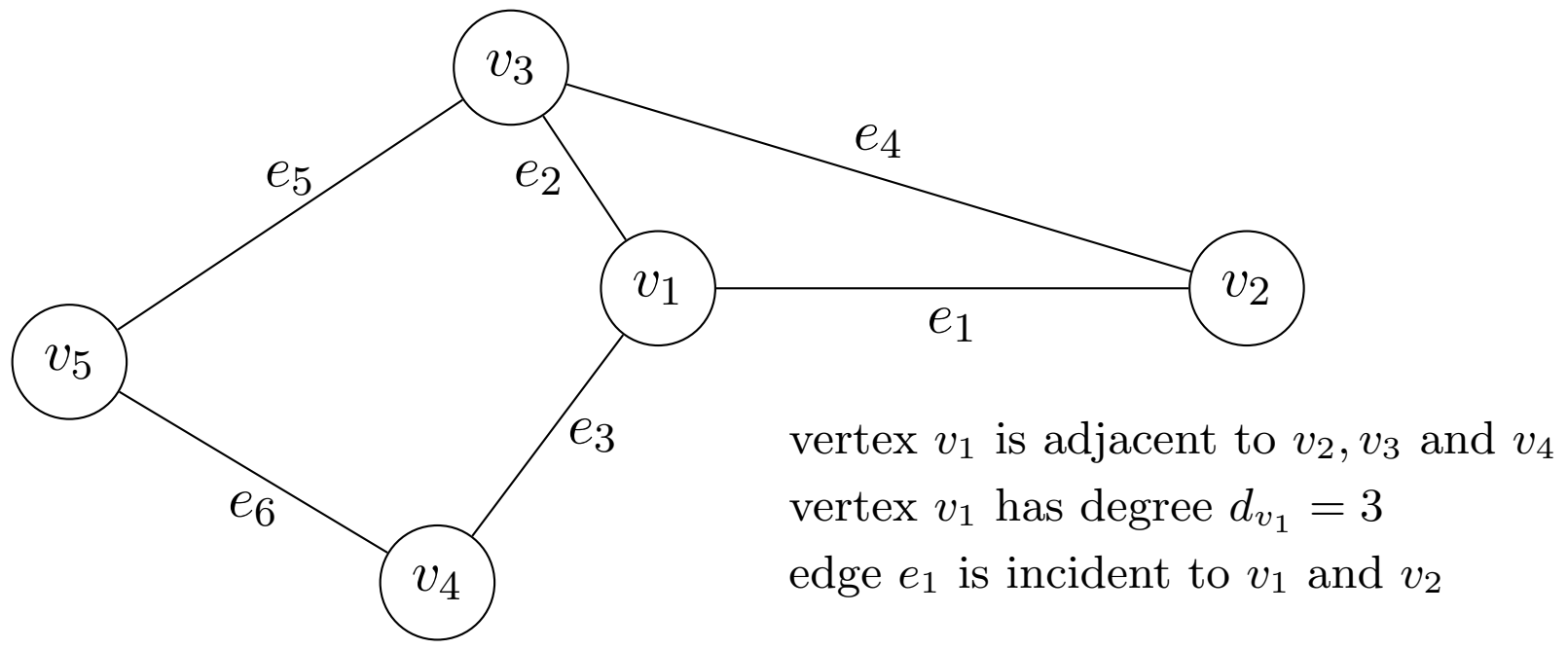

A network or \textit{graph} is a pair $G=(V,E),$ where $V = \lbrace v_1, \ldots, v_n \rbrace$ is a set of \textit{vertices} and $E= \lbrace e_1, \ldots, e_m \rbrace$ a set of paired vertices (\textit{edges}). A first example of a graph is given in Figure \ref{fig:ex-graph}. Here, $G$ has five vertices $v_1, \ldots, v_5$ and six edges $e_1, \ldots, e_6.$ Two vertices are \textit{adjacent} if they are connected by an edge. These vertices are also called \textit{neighbors}. An edge is \textit{incident} to a vertex $v$ if the edge is
attached to $v$.

Besides the number of vertices $n$ and the number of edges $m$, an important characteristic of a graph is the degree. The \textit{degree} $d_v$ of a vertex $v$ is defined as the number of neighbors of a vertex
which is the same as the number of incident edges.
 
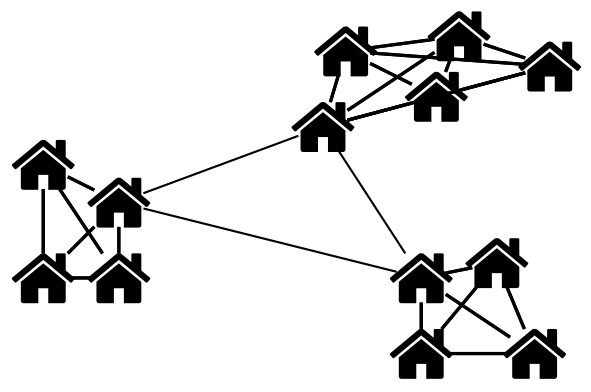
\begin{figure}[ht]
\begin{center}
\begin{tikzpicture}[scale=0.50]
\node[scale=1.2] (A1) at (3.5,3.2) {};
\node[scale=1.2] (A2) at (1.7,2) {};
\node[scale=1.2] (A3) at (3.2,2.4) {};
\node[scale=1.2] (A4) at (4.7,2.8) {};
\node[scale=1.2] (A5) at (2,3) {};
\node[scale=1.2] (B1) at (-1,0) {};
\node[scale=1.2] (B2) at (-2,0) {};
\node[scale=1.2] (B3) at (-1,1) {};
\node[scale=1.2] (B4) at (-2,1.5) {};
\node[scale=1.2] (C1) at (2+2.5,-2.5+1.5) {};
\node[scale=1.2] (C2) at (2+1,-2.5+1.5) {};
\node[scale=1.2] (C3) at (1+2,-1+1) {};
\node[scale=1.2] (C4) at (2+2,-1.3+1.5) {};
\foreach \i in {1,...,5}
{
    \foreach \j in {1,...,5}
    {
        \draw[semithick] (A\i) -- (A\j);   
    }
}
\foreach \i in {1,...,4}
{
    \foreach \j in {1,...,4}
    {
        \draw[semithick] (B\i) -- (B\j);   
        \draw[semithick] (C\i) -- (C\j);   
    }
}
\draw[] (B3) -- (C3);
\draw[] (B3) -- (A2);
\draw[] (C3) -- (A2);
\node[scale=1.2] at (3.5,3.2) {\faHome};
\node[scale=1.2] at (1.7,2) {\faHome};
\node[scale=1.2] at (3.2,2.4) {\faHome};
\node[scale=1.2] at (4.7,2.8) {\faHome};
\node[scale=1.2] at (2,3) {\faHome};
\node[scale=1.2] at (-1,0) {\faHome};
\node[scale=1.2] at (-2,0) {\faHome};
\node[scale=1.2] at (-1,1) {\faHome};
\node[scale=1.2] at (-2,1.5) {\faHome};
\node[scale=1.2] at (2+2.5,-2.5+1.5) {\faHome};
\node[scale=1.2] at (2+1,-2.5+1.5) {\faHome};
\node[scale=1.2] at (1+2,-1+1) {\faHome};
\node[scale=1.2] at (2+2,-1.3+1.5) {\faHome};

\end{tikzpicture}
\caption{Houses and communities.}
\label{fig:ex-graph-village}
\end{center}
\end{figure}

Another example: Imagine the vertices of a graph to be houses on an island. Several houses are
connected by streets, hence, the streets are the edges of our graph. Along the streets, houses
can exchange things with each other, such as agricultural products, information, etc. Like
in Figure \ref{fig:ex-graph-village}, we might identify groups of houses on the island that are
more densely connected to each other than to other houses; we could speak of villages. In graph
theoretical terms, such densely interconnected sets of vertices are called \textit{communities}.

Whether a network is organized in community structures can be investigated, for example, by the Louvain community finding algorithm. 
This algorithm can also calculate the \textit{modularity} of a graph, a measure to which degree such an organization is present. 

\section{Brain Network Model}

The first idea for our model is to represent the human brain as a network consisting of vertices
and edges with neurons as vertices and synapses as edges.
Unfortunately, the resolution of common MRI scanners is not high enough to render the extremely small
neurons. Instead, a common approach is to consider groups of neurons that are close to each other as
\textit{brain regions}.

\begin{figure}[ht]
     \centering
     \begin{subfigure}[b]{0.6\textwidth}
         \centering
         \includegraphics[width=\textwidth]{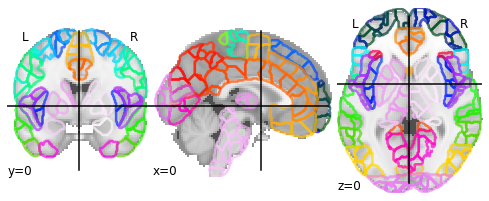}
         \caption{Brain atlas with 571 brain regions}
         \label{fig:atlas}
     \end{subfigure}
     \\
     \begin{subfigure}[b]{0.6\textwidth}
         \centering
         \includegraphics[width=\textwidth]{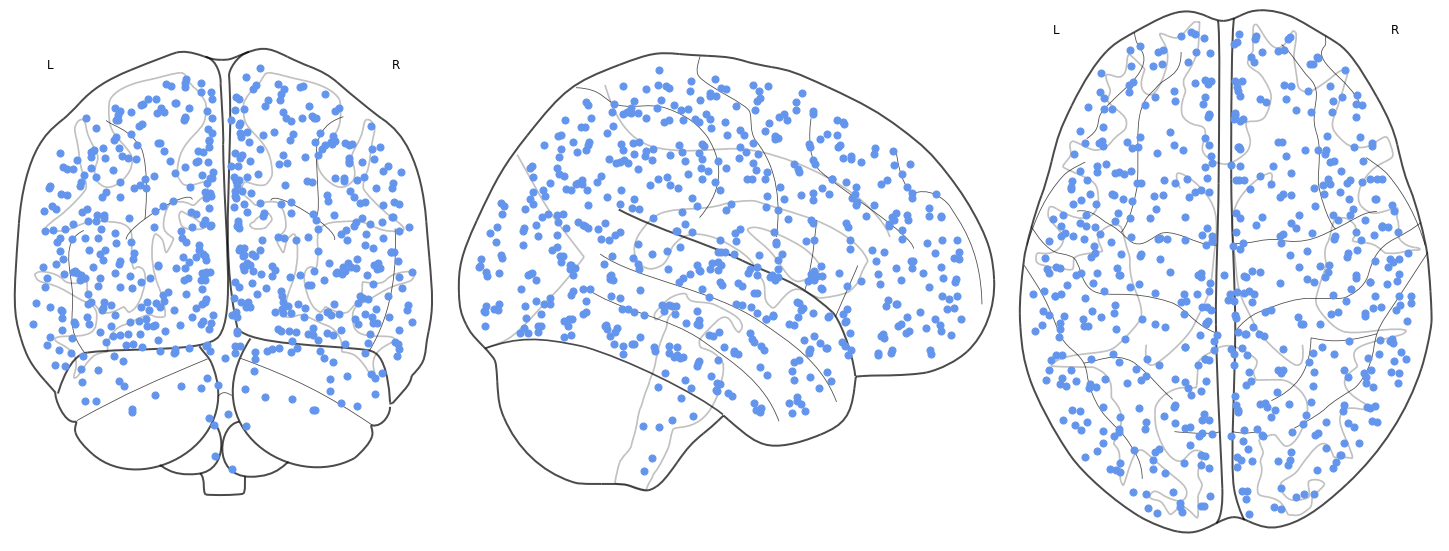}
         \caption{Vertices of the brain network}
         \label{fig:vertices}
     \end{subfigure}
     \\
     \begin{subfigure}[b]{0.6\textwidth}
         \centering
         \includegraphics[width=\textwidth]{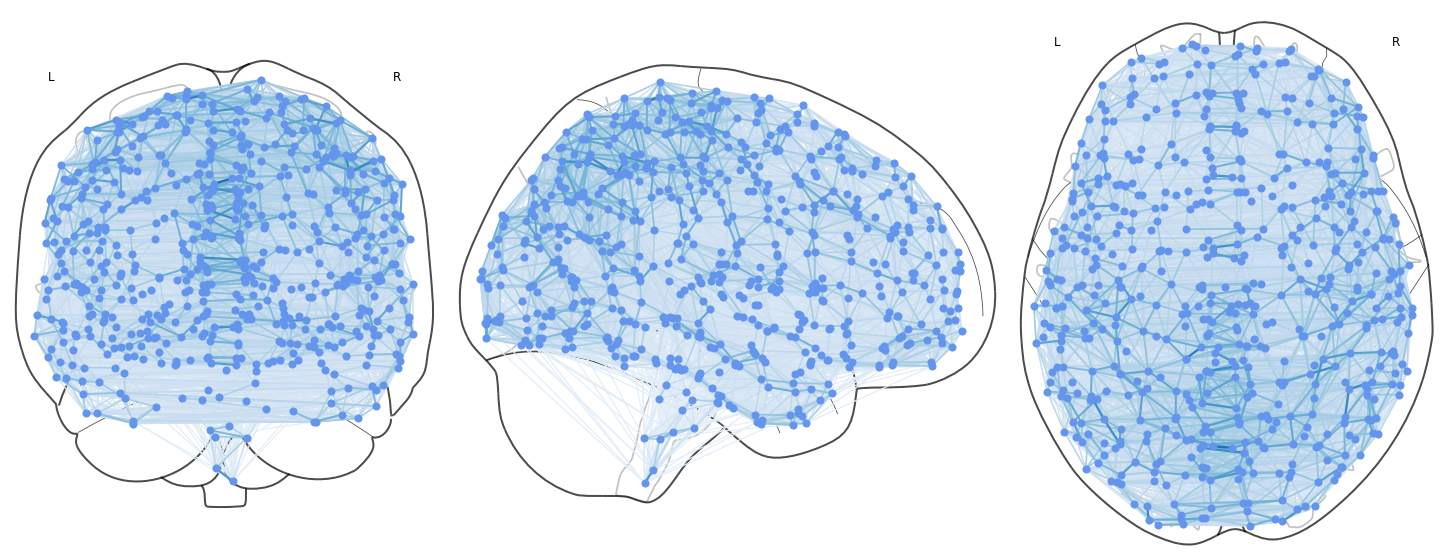}
         \caption{Vertices and edges of the brain network}
         \label{fig:edges}
     \end{subfigure}
        \caption{Construction of the Brain Network Model.}
        \label{fig:brain-network-model}
\end{figure}

For the brain as a three-dimensional object, a definition of these brain regions may be achieved by a partitioning into adjacent regions that have approximately the same size and shape. 
In our model, this resulted in the parcellation displayed in Figure \ref{fig:atlas} with 571 regions with an approx.\ volume of 2 ml \cite{WB}. 
Each of these regions now serves as a vertex in our graph, where we position the vertex in the middle of the brain region, see Figure \ref{fig:vertices}.

The resulting groups of neurons can be structurally or functionally connected. A structural connection is a physical connection through fibres in the brain. A functional connection is more likely to refer to the collaboration of different regions. If two areas in the brain are often active at the same time, for instance, to perform a task, they are said to be functionally connected. 

Depending on which connectivity measure one considers, one examines the structural or functional brain network. 
In both cases, the brain regions form the vertex set of our graph. 
The edges, or streets, connecting the vertices are either the structural or functional connections (see Figure \ref{fig:edges} for an example of a functional connectome).

Modeling the brain as a graph enables us to apply new methods for analyzing its structure and
organizing principles. A variety of algorithms from graph theory can now be applied to identify
regions with particularly strong connections, to determine shortest paths, or to find communities.
Several authors have performed graph theoretical analysis of the brain network. For example,
it was found that healthy brain networks are generally organized, respecting economic principles
of minimizing wiring costs and maximizing global efficiency, see, for example, \cite{BS12} and the
references therein. 

If every two houses on our imaginary island were to be connected by a street, this would cause
massive construction costs. An economic trade-off would be a high connectivity in local
communities (villages), where wiring costs are comparably low, combined with only a few expensive
long-distance roads (expressways) to more remote communities.
This topological property is characterized by high modularity and is a particular kind of a
\textit{small-world} network architecture allowing highly efficient information transfer with
low wiring costs \cite{BS09}.

By contrast, in many neurological disorders, peculiarities in the brain network have been noticed,
for example, weakened connectivity or modified organizational structures.
In \cite{WB}, the community formation in brain networks of a group of healthy adults were compared to a group of AD patients. 
The results of this study suggest that the community structures in the Alzheimer's group have a significantly lower modularity than the healthy control group. 
This decrease of modularity in AD patients compared with both young and senior controls indicates a loss of efficiency of information transfer across the whole brain.

\section{Spreading of Tau Proteins}

The second important part of the Computational Brainsphere Model is the Tau Distribution Model which aims to simulate the spreading of tau proteins across the brain network. 
To simulate such processes, we use partial differential equations (PDEs). 
A partial differential equation combines partial derivatives of a function in several dimensions, and the aim is to find this function from this equation.
In our first mathematical model, these two dimensions are time and space. We want to simulate how the concentration of tau at position $x$ in the network changes over time $t$. 
Following the transneuronal spread hypothesis, we suppose that tau proteins
can only spread along the edges of the brain network. We will here
model this spreading as a diffusion process.

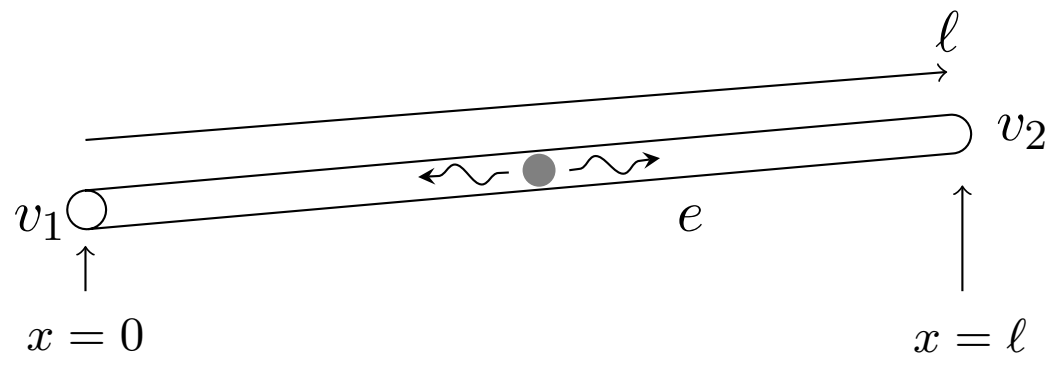
\begin{figure}
\begin{center}
\begin{tikzpicture}
\node[cylinder, draw, shape aspect=1,rotate=185,minimum height=6cm,fill=white,fill opacity=1]{} ;
\draw[black] (1,-0.1) node[anchor=north]{$e$};
\draw[black] (-3.3,-0.1) node[anchor=north]{$v_1$};
\draw[black] (3.2,0.5) node[anchor=north]{$v_2$};
\draw[->] (-3,0.2) -- (2.7,0.65) node[anchor=south]{$\ell$};
\draw[->] (-3,-0.8) -- (-3,-0.5);
\draw[black] (-3,-0.85) node[anchor=north]{\footnotesize $x=0$};
\draw[->] (2.8,-0.8) -- (2.8,-0.1);
\draw[black] (2.85,-0.85) node[anchor=north]{\footnotesize $x=\ell$};
\filldraw[gray] (0,0) circle (2.9pt);
\draw[-stealth,decorate,decoration={snake,amplitude=2pt,pre length=2pt,post length=3pt}] (-0.2,-0.015) -- ++ (-0.6,-0.03);
\draw[-stealth,decorate,decoration={snake,amplitude=2pt,pre length=2pt,post length=3pt}] (0.2,0) -- ++ (0.6,0.08);
\end{tikzpicture}
\end{center}
\caption{Diffusion across a single, isolated edge.}
\label{fig:single-edge}
\end{figure}

To explain this, consider first only a single, isolated edge of the network, for example, the edge $e=(v_1,v_2)$ between vertex $v_1$ and vertex $v_2$. 
Physically, we imagine the edge to be a thin tube with length $\ell$, see Figure \ref{fig:single-edge}.
Vertex $v_1$ is at position $x=0$ and $v_2$ at position $x= \ell$ of the edge. The concentration of tau at any position $x$ of the edge at time $t$ is given by a function denoted by 
$\tau =\tau(x,t)$ 
which describes the amount of tau per unit volume. 

Fick's first law tells us that the concentration $\tau$ can be related to the flux $F$ by the equation
$$
F = -D \cdot \frac{\partial \tau}{\partial x},
$$
where $D$ is a (given) diffusion coefficient. 
By flux, we mean the amount of a substance per unit area per unit time. 
The expression $\frac{\partial \tau}{\partial x}$ stands for the first partial derivative of $\tau$ with respect to $x$, i.e., it describes how the concentration changes with respect to the space
coordinate.

Using this relation and respecting mass conservation, we can derive an equation that tells us how diffusion causes the concentration to change with respect to time. This equation is Fick's second law and reads
\begin{equation}\label{Eq:Fick-2nd-law}
\frac{ \partial \tau}{\partial t} = D \cdot \frac{\partial^2 \tau}{\partial x^2}.
\end{equation}
Here  $\frac{ \partial \tau}{\partial t}$ denotes the first partial derivative of $\tau$
with respect to time. Equation (1) is a partial differential equation (PDE) as it involves
derivatives of $\tau$ with respect to both time and space. Solving this PDEs
tells us how the concentration of $\tau$ at any position $x$ on an edge changes with respect
to time $t$.

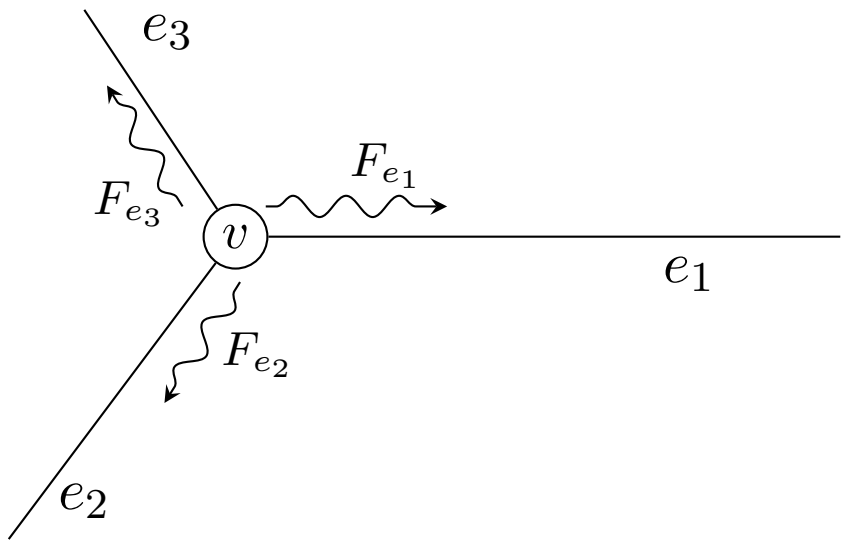
\begin{figure}
\begin{center}
\begin{tikzpicture}
\draw[black] (3,0) node[anchor=north]{$e_{1}$};
\draw[black] (0,0) -- (4,0);
\draw[-stealth,decorate,decoration={snake,amplitude=2pt,pre length=2pt,post length=3pt}] (0.2,0.2) -- ++(1.2,0);
\draw[black] (1,0.75) node[anchor=north]{\footnotesize $F_{e_1}$};
\draw[black] (-1,-1.5) node[anchor=north]{$e_{2}$};
\draw[black] (-1.5,-2) -- (0,0);
\draw[-stealth,decorate,decoration={snake,amplitude=2pt,pre length=2pt,post length=3pt}] (0.03,-0.3) -- ++(-0.5,-0.8);
\draw[black] (0.15,-0.5) node[anchor=north]{\footnotesize $F_{e_2}$};
\draw[black] (-0.45,1.6) node[anchor=north]{$e_{3}$};
\draw[black] (-1,1.5) -- (0,0);
\draw[-stealth,decorate,decoration={snake,amplitude=2pt,pre length=2pt,post length=3pt}] (-0.35,0.2) -- ++(-0.5,0.8);
\draw[black] (-0.7,0.5) node[anchor=north]{\footnotesize $F_{e_3}$};
\filldraw[white] (0,0) circle (6pt);
\draw[black] (0,0) circle (6pt);
\draw[black] (0,0) node{\small $v$};
\end{tikzpicture}
\end{center}
\caption{Diffusion at a vertex $v$ with three incident edges $e_1,e_2,e_3$}
\label{fig:diffusion_vertex}
\end{figure}

Now that we elaborated how we can describe the diffusion process on a single edge, we can attend to the generalization to the whole network. 
Consider Figure \ref{fig:diffusion_vertex}, where we illustrate as an example one vertex $v$ with three incident edges $e_1, e_2$ and $e_3$.
In general, a vertex $v$ has $d_v$ incident edges $e_1, \ldots, e_{d_v}$ as explained in section
\ref{sec:graph-theory}. 
For each of these edges, we now need a separate function to describe the concentration along the coordinates of the edge. 
We will, therefore, denote the function for edge $e_i$ by $\tau_i$.

For ease of notation, we suppose that the vertex $v$ lies at position $x=0$ for every attached edge $e_i$. 
This means that the concentration of tau at vertex $v$ is given by $\tau(v) := \tau_i(0)$. 
This definition only makes sense if $\tau_1(0) = \tau_2(0) = \ldots = \tau_d(0)$ for all incident edges which we assume in the following (\textit{continuity condition}).

On every edge attached to vertex $v$, the diffusion can be described by equation (\ref{Eq:Fick-2nd-law}), but what happens at the vertex $v$ itself? According to Fick's first law, the rate at which tau leaves node $v$ along an edge $e_{i}$ can be modelled as
\begin{equation}\label{flow_edge}
F_{e_i}(t)=- D \cdot \frac{\partial \tau_{i}}{\partial x} \Big|_{x=0},
\end{equation}
where we evaluate the derivative at the point $x=0$ since we agreed that the vertex is at position $x=0$ on every attached edge.  Hence, the total flux out of node $v$ along all incident edges $e_1, \ldots, e_{d_v}$ is given by the sum of the flows over the single edges, i.e.,
\begin{equation*}\label{flow_node}
F_{v}(t)= -D \cdot \sum_{i = 1}^{d_v} \frac{\partial \tau_{i}}{\partial x} \Big|_{x=0}.
\end{equation*}

If we assume that no substance is produced or lost in the vertices, the total flux out of vertex $v$ must be zero as everything that flows in must flow out again. This is also known as conservation of currents and provides us with the so-called Neumann-Kirchhoff vertex condition. It couples the diffusion equations
on the different edges at the vertices \cite{BK}.
Together, we arrive at a coupled system of partial differential equations that describes a diffusion process on a network. This system is given by
\begin{equation}\label{eq:tau_diffusion}
\begin{aligned}
& \hspace*{0.1cm} \frac{ \partial \tau_i}{\partial t} =  D\, \frac{\partial^2 \tau_i}{\partial x^2} \hspace*{1.0cm} \text{for all edges } e_i = e_1, \ldots, e_m  \\
& \left. \begin{aligned}
        &\tau(v) = \tau_i(0) \\
        &\displaystyle \sum_{i = 1}^{d_v} \frac{\partial \tau_{i}}{\partial x} \Big|_{x=0} = 0  
\end{aligned}  \right\} \text{ for all vertices } v \text{ with incident edges } e_i = e_1, \ldots, e_{d_v}.
\end{aligned}
\end{equation}
%
\section{Numerical Simulations}
The goal of a numerical simulation is to predict patterns of future tau accumulation, given an initial tau pattern that was actually measured in the brain of a patient.   
In other words, given the initial tau distribution at time $t=0$, we want to solve system \eqref{eq:tau_diffusion} to obtain the tau concentration $\tau_i(x,t)$ for a future time point $t > 0$ and at each position $x$. This can be achieved by numerical methods redesigned for graphs, such as finite differences, finite element or spectral methods, see, for example, \cite{AW,AB,W}.

\begin{figure}[h]
	\centering
	\includegraphics[width=0.8\textwidth]{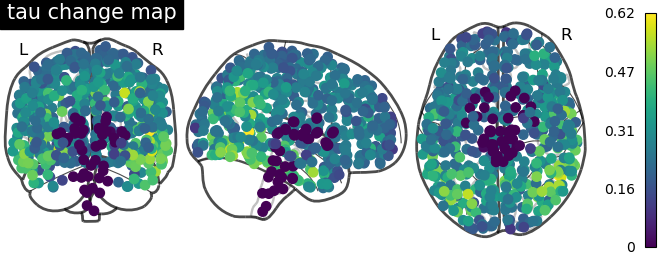}
	\includegraphics[width=0.8\textwidth]{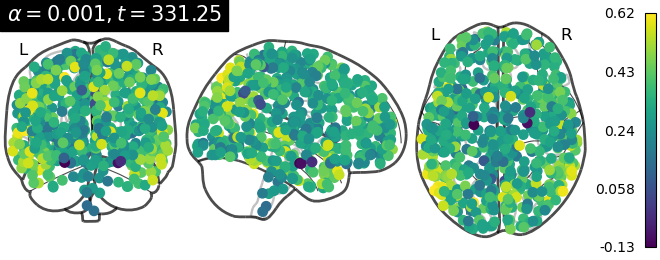}
	\includegraphics[width=0.8\textwidth]{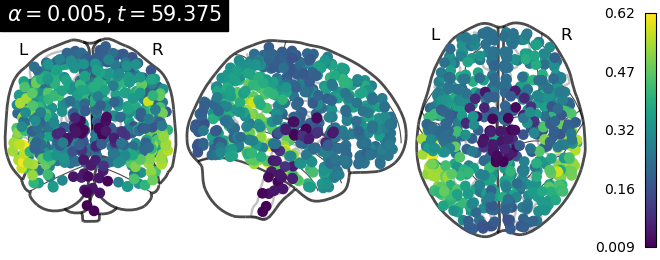}
	\caption{Actual tau change map and deviation of the predicted tau pattern and the initial data for $\alpha=0.001,0.005$. \textit{Graphical representation using \texttt{nilearn} in python} \cite{W}.}\label{fig:tau_simulation}
\end{figure}

In fact, \eqref{eq:tau_diffusion} is only a very simple attempt to describe the spread of tau proteins across the brain network. 
In a more realistic model, other processes such as chemical reactions with the network environment are important. 
We therefore now present results for the solution of a so-called reaction-diffusion equation
$$
\frac{ \partial \tau_i}{\partial t} =  D\, \frac{\partial^2 \tau_i}{\partial x^2} + \alpha \tau_i \quad \hspace*{1.0cm} \text{for all edges } e_i = e_1, \ldots, e_m 
$$
where $\alpha$ is some given constant. 
We compute the deviation of the predicted tau pattern with respect to the initial condition in every brain region and compare the predicted change to the actual tau change measured approximately one year after the initial scan. 
The predictions in Figure \ref{fig:tau_simulation} were obtained using a finite element discretization; the exact discretization details and a description of the selected data can be found in \cite{W}. 
The results suggest that the reaction term plays an important role in the prediction of tau patterns.

\end{document}